\documentclass[letterpaper]{amsart}
\usepackage{amsmath}
\usepackage{amsthm}
\usepackage{amssymb}
\usepackage{amscd}          			%diagrams
\usepackage{color}  
\usepackage{bbm}            			%bold numbers
\usepackage{mathrsfs}           		%nice calligraphy
\usepackage[bookmarksnumbered,plainpages,backref]{hyperref} % Laurent i added new package

\usepackage[french,english]{babel}     		%if writing in French

\usepackage{graphicx}          		%uncomment for inclusion of graphic images

\usepackage{verbatim}           		%DO NOT comment this line

\usepackage{MnSymbol}

\theoremstyle{plain}
\newtheorem{Theorem}{Theorem}[section]
\newtheorem{Lemma}[Theorem]{Lemma}
\newtheorem{Corollary}[Theorem]{Corollary}
\newtheorem{Proposition}[Theorem]{Proposition}
\newtheorem{Conjecture}[Theorem]{Conjecture}
\newtheorem{theorem}[Theorem]{Theorem}
\newtheorem{proposition}[Theorem]{Proposition}

\newtheorem{lemma}[Theorem]{Lemma}

\theoremstyle{definition}
\newtheorem{Definition}[Theorem]{Definition}
\newtheorem{Example}[Theorem]{Example}
\newtheorem{definition}[Theorem]{Definition}

\theoremstyle{remark}
\newtheorem{Remark}[Theorem]{Remark}
\newtheorem{Claim}[Theorem]{Claim}

\newcommand{\bbC}{\mathbb{C}}       %complex numbers
\newcommand{\N}{\mathbb{N}}     %natural numbers
\newcommand{\R}{\mathbb{R}}     %real numbers
\newcommand{\Z}{\mathbb{Z}}         %integers

\newcommand{\calD}{\mathscr{D}}

\DeclareMathOperator{\dist}{dist}			%distance
\DeclareMathOperator{\diver}{div}			%divergence
\DeclareMathOperator{\im}{im}				%image
\DeclareMathOperator{\supp}{supp}			%support
\DeclareMathOperator{\spann}{span}			%support
\usepackage[mathcal]{eucal}					%for the set of BV subsets of Rn

\newcommand{\bsC}{{\mathcal{C}}}

\newcommand{\bsL}{{\mathcal{L}}}

\newcommand{\bsT}{{\mathcal{T}}}

\newcommand{\LL}{\mathrm{L}}

\newcommand{\eps}{\epsilon}

\renewcommand{\geq}{\geqslant}
\renewcommand{\leq}{\leqslant}
\renewcommand{\ge}{\geqslant}
\renewcommand{\le}{\leqslant}
\renewcommand{\epsilon}{\varepsilon}

\newcommand{\liminfe}{\mathop{\underline{\lim}}}
\newcommand{\limsupe}{\mathop{\overline{\lim}}}

\begin{document}

\title[Divergence-type equations associated to elliptic systems]{Continuous solutions for divergence-type equations associated to elliptic systems of complex vector fields}
\author{Laurent Moonens and Tiago Picon}
\date{\today}
\thanks{Laurent Moonens and Tiago Picon were partially supported by the French ANR project ``GEOMETRYA'' no.~ANR-12-BS01-0014 and S\~ao Paulo Research Fundation - Fapesp grant 2013/17636-5, respectively.}
%\begin{comment}

\maketitle
\begin{abstract}
{In this paper, we characterize all the distributions $F \in \calD'(U)$ such that there exists a continuous weak solution $v \in C(U,\mathbb{C}^{n})$ (with $U \subset \Omega$) to the divergence-type equation $$L_{1}^{*}v_{1}+...+L_{n}^{*}v_{n}=F,$$
where  $\left\{L_{1},\dots,L_{n}\right\}$ is an elliptic system of linearly independent vector fields with smooth complex coefficients defined on $\Omega \subset \R^{N}$. In case where $(L_1,\dots, L_n)$ is the usual gradient field on $\R^N$, we recover the classical result for the divergence equation proved by T. De Pauw and W. Pfeffer.}

\end{abstract}
\bigskip

%Suppose that $\bsL=\left\{L_{1},\dots,L_{n}\right\}$ is a  system of linearly independent vector fields with smooth complex coefficients defined on an open set $\Omega \subset \R^{N}$. We may consider the operators $\nabla_{\bsL}\,u\doteq(L_{1}u,\dots,L_{n}u)$ for $u\in C^{\infty}(\Omega)$ and ${\rm div}_{\bsL}\,v\doteq\sum_{j=1}^{n}L_{j}v_{j}$ for $v\in C^{\infty}(\Omega,\Lambda^{1}\R^{n})$  which are precisely the operators $\nabla$ and $\rm{div}$ when $n=N$ and $L_{j}=\partial_{x_{j}}$. 

\section{Introduction}

Recently a series of new results on the classical divergence equation have been published. In the original paper due to J. Bourgain and H. Brezis \cite{BB1} the authors presented new developments for the solvability of the equation
\begin{equation}\label{divergente}
\diver v =F,
\end{equation}
%on the torus $\mathbb{T}^{N}$ for $N \geq 2$ concerning solvability results 
when $F \in L_{\#}^{p}(\mathbb{T}^{N})=\left\{ f \in L^{p}(\mathbb{T}^{N})\; |\; \int_{\mathbb{T}^{N}}f=0 \right\}$, in the special limiting case $p=N$. A surprising result \cite[Theorem~1']{BB1} asserts that for every $f \in L^{N}_{\#}(\mathbb{T}^{N})$ there exists a continuous solution of \eqref{divergente}.

Concerning continuous solutions to (\ref{divergente}) in the whole Euclidean space, T. de Pauw and W. Pfeffer \cite{PP} characterized the (real) distributions $F$ for which the equation (\ref{divergente}) has a continuous solution, \emph{i.e.}  there exists $v\in C(\R^N,\R^N)$ such that the following holds:
$$F(\varphi)=-\int_{\R^{N}}v\cdot \nabla \varphi,$$ %, \quad \forall\; \varphi \in \calD(\R^{N}), $$
for every test function $\varphi \in \calD(\R^{N})$. They show such distributions are exactly the ones satisfying a particular continuity property: for each $\epsilon > 0$ there should exist a constant $\theta>0$ such that one has:
\begin{equation}\label{strongcharge}
|F(\varphi)|\leq \theta\|\varphi\|_{{1}}+\epsilon \|\nabla \varphi \|_{{1}},
\end{equation}    
for all $\varphi\in\calD(\R^N)$ supported in the ball centered at the origin with radius $1/\epsilon$. As a particular case, they show that (the distribution associated to) any function $F\in L^N(\R^N)$ enjoys that property, so that in particular (\ref{divergente}) is continuously solvable for all $F\in L^N(\R^N)$.

{Integral estimates in $L^{1}$ norm like  \eqref{strongcharge} have been studied in several settings, among which div-curl and elliptic-canceling operators, measure and divergence-free vector fields, nilpotent groups, CR complexes and applications to fluid dynamics. We refer to \cite{VS5} for an overview and development of these subjects.}

{The results obtained previously for \eqref{divergente}} are closely related to the gradient $\nabla$ generated by the canonical vector fields {$L_{j}=\partial_{x_j}$} for $j=1,...,N$. Suppose now that $\bsL:=\left\{L_{1},\dots,L_{n}\right\}$ is a system of linearly independent vector fields with smooth complex coefficients defined on an open set $\Omega \subset \R^{N}$. {Analogously,} we may consider the gradient associated to the {system $\bsL$} defined by  
$\nabla_{\bsL}\,u :=(L_{1}u,\dots,L_{n}u),$ for $u\in C^{\infty}(\Omega)$ and its formal {complex}  adjoint operator 
\begin{equation}\label{diverl}
{\rm div}_{\bsL^{*}}\,v:=\sum_{j=1}^{n}L^{*}_{j}v_{j}, \quad v\in C^{\infty}(\Omega,\mathbb{C}^{n}),
\end{equation}
 which are precisely the operators $\nabla$ and $\diver$ when $n=N$ and $L_{j}=\partial_{x_{j}}$. We use the notation  $L_{j}^{*}\;:=\;\overline{ L_j^{t}}$ where $\overline L_{j}$ denotes the vector field obtained from $L_j$ by conjugating its coefficients and $L_j^t$ is the formal transpose of $L_j$ for $j=1,\dots,n$~---~namely this means that, for all (complex valued) $\varphi,\psi\in\calD(\Omega)$, we have:
 $$
 \int_\Omega (L_j\varphi)\bar{\psi}=\int_\Omega \varphi \overline{L_j^*\psi}.
 $$
 
The following version of the $L^{1}$ Sobolev-Gagliardo-Nirenberg theorem associated to $\nabla_{\bsL}$ was proved in \cite{HP1}, namely:
\begin{theorem}\label{thmHP}
Assume that the system of vector fields $L_{1},...,L_{n}$, $n \geq 2$, is linearly independent  and elliptic. Then every point $x_{0} \in \Omega$ is contained in an open neighborhood $U \subset \Omega$ such that
\begin{equation}\label{SGN}
\|\varphi\|_{L^{N/N-1}}\leq C\sum_{j=1}^{n}\|L_{j}\varphi\|_{L^{1}}, \;\;\;\forall\;\varphi \in \calD (U),
\end{equation}
holds for $C=C(U)>0$. Conversely, if \eqref{SGN} holds then the system must be elliptic on $U$. 
\end{theorem}
%\textcolor{blue}{I put this sentence in the Section 2: Here ellipticity means that for any \textit{real} $1$-form $\omega$  such that  $\langle\omega,L_j\rangle=0$ holds for $j=1,..,n$, we have $\omega=0$; equivalently, this also means that the second-order operator $\Delta_\bsL=\diver_{\bsL^*}(\nabla_\bsL \cdot)$ is elliptic.}

In this work we are interested to study the (local) continuous solvability of the equation:
\begin{equation}\label{diverl}
\diver_{\bsL^{*}} v =F.
\end{equation}
Our main result is the following.
\begin{theorem}\label{mainthm}
Assume that the system of vector fields $L_{1},...,L_{n}$, $n \geq 2$, is linearly independent  and elliptic. Then every point $x_{0} \in \Omega$ is contained in an open neighborhood $U \subset \Omega$ such that for any $F \in \calD'({U})$, the equation (\ref{diverl}) is continuously solvable in $U$ if and only if $F$ is an $\bsL$-charge in $U$, meaning that for every $\epsilon>0$ and every compact set $K \subset \subset U$, there exists $\theta=\theta(K,\epsilon)>0$ such that one has:%for every  i$\varphi \in \calD'{K}$,
\begin{equation}\label{lstrong}
\left| F(\varphi) \right| \leq C\|\varphi\|_{{1}}+\epsilon\|\nabla_{\bsL} \varphi\|_{{1}}, 
\end{equation}
for all $\varphi\in\calD_K(U)$~---~the latter being the set of all smooth functions in $U$ supported inside $K$.
\end{theorem}

One simple argument (see Section \ref{sec2}) shows that the above continuity property on $F$ is a necessary condition for the continuous solvability of equation \eqref{diverl} in $U$. Theorem \ref{mainthm} asserts that the continuity property \eqref{lstrong} is also sufficient, under the ellipticity assumption on the system of vector fields.

The organization of the paper is as follows. In Section \ref{sec0}, we study some properties of elliptic systems of complex vector fields. Section \ref{sec1} is devoted to the definition and some properties of the space $BV_{\bsL,c}$ of functions with bounded  $\bsL$-variation. In Section \ref{sec2}, we discuss linear functionals on $BV_{\bsL,c}$ called $\bsL$-charges. The proof of our main result is presented in Section \ref{sec3}. The Appendix is concerned with technical results on   pseudodifferential operators, mainly on their boundedness and compactness.    

\vglue 0.5cm

{\textbf{Notations.} We always denote by $\Omega$ an open set of $\R^{N}$, $N\geq 2$. Unless otherwise specified, all functions are complex valued and the notation $\int_A f$ stands for the Lebesgue integral $\int_A f(x)dx$.  As usual, $\calD(\Omega)$ and $\calD'(\Omega)$ are the spaces of complex test functions and distributions, respectively. When $K\subset \subset \Omega$ is a compact subset of $\Omega$, we let $\calD_{K}(\Omega):=\calD(\Omega)\cap\mathcal{E}'(K) $, where $\mathcal{E}'(K)$ is the space of all distributions with compact support in $K$. Since the ambient field is $\bbC$, we identify (formally) each $f\in L^1_{loc}(\Omega)$  with the distribution $T_f\in\calD'(\Omega)$ given by $T_f(\varphi)=\int_\Omega f\bar{\varphi}$. We consider $C(\Omega,\mathbb{C}^{n})$ the space of all continuous vector-valued functions $v:\Omega \rightarrow \mathbb{C}^{n}$. We also introduce the notation $\|\nabla_{\bsL}\varphi\|_{p}:=\sum_{j=1}^{n}\|L_{j}\varphi\|_{p}$ (where $\|\cdot\|_{p}$ is the standard norm in $L^{p}(\Omega)$) for $1\leq p\leq \infty$. Finally we use the notation $f\lesssim g$ to indicate the existence of an universal constant $C>0$, independent of all variables and unmentioned parameters, such that one has $f \leq Cg$.}

%We now study some properties of elliptic systems of vector fields.

%
%
\section{Ellipticity and its consequences}\label{sec0}

Consider $n$ complex vector fields $L_1,\dots,L_n$, $n\ge1$, with smooth coefficients defined on a neighborhood $\Omega$ of the origin in $\R^N$, $N\ge2$. We will assume that the vector fields $L_{1},...,L_{n}$ do not vanish in $\Omega$, in particular, they may be viewed as nonvanishing sections of the vector bundles  $\mathbb{C} T\Omega$  as well as first order differential operators of principal type. %\textcolor{blue}{sugestion to remove:The Sobolev-Gagliardo-Nirenberg estimate \eqref{SGN} above expresses a property of this bundle rather than a property of this specific set of generators.}

In the sequel, we will always assume (unless otherwise mentioned) that the following two properties hold:
\begin{itemize}
  \item[(a)] $L_1,\dots,L_n$ are everywhere linearly independent;
  \item[(b)] the system $\{L_1,\dots,L_n\}$ is \textit{elliptic}.
\end{itemize} 
The latter means for any 1-form $\omega$ (i.e. any section of $T^{*}(\Omega)$), the equality $\langle \omega,L_{j} \rangle=0$ for $1 \leq j \leq n$ implies that one has $\omega=0$. Consequently, the number $n$ of vector fields must satisfy $\frac{N}{2}\le n\le N$   
\footnote {In fact, if one writes $L_{j}=X_{j}+iY_{j}$ where $\left\{X_{j},Y_{j} \right\}_{j}$ are real vector  fields, then $2n \geq N$. Suppose indeed that $\#\left\{X_{j},Y_{j} \right\}=2n<N$. Then there exist $f \notin \spann \left\{X_{j},Y_{j} \right\}_{j}$ and $\omega :=d_{f}\neq 0$ such that $\omega(L_{j})=0$ for $j=1,...,n$ but $\omega \neq 0$; that is a contradiction, since  the system $\{L_1,\dots,L_n\}$ is supposed to be \textit{elliptic}. Clearly, on the other hand, we have $n\leq N$.} . 
Alternatively the assumption (b) is equivalent to require that the second order operator 
\begin{equation}\label{deltal}
\Delta_{\bsL}:=L_{1}^{*}L_{1}+...+L_{n}^{*}L_{n}=\diver_{\bsL^{*}}\nabla_{\bsL}
\end{equation}
is elliptic. Using a representation of vector fields in local coordinates $(x_{1},...,x_{N})$ we can  assume that one has:
\begin{equation}
L_{j}=\sum_{k=1}^{N}c_{jk}\partial_{x_{k}}\quad j=1,...,n,
\end{equation}
with smooth coefficients globally defined on $\R^{N}$ that possess bounded derivatives of all orders. A simple computation implies then that one has $L^{*}_{j}=-\overline{L_{j}}+c_{j}$ where $c_{j}:=\sum_{k=1}^{N}\partial_{x_{k}}c_{jk}$; the (uniform) ellipticity means that there exists $c>0$ such that one has
$$\sum_{j=1}^{n}\left| \sum_{k=1}^{N}c_{jk}(x)\xi_{k} \right|^{2} \geq c |\xi|^{2}, $$ 
for all $x,\xi \in \R^{N}$.

The second-order (elliptic) operator $\Delta_{\bsL}$ may be regarded as an elliptic pseudodifferential operator with symbol in the H\"ormander class $S^{2}_{1,0}(\Omega)$. Hence there exist {scalar-valued} properly supported pseudodifferential operators $q(x,D) \in OpS^{-2}_{1,0}(\Omega)$ and $r(x, D) \in OpS^{-\infty}(\Omega)$ such that one has:
\begin{equation}\label{pseudo}
\Delta_{L}q(x,D)f+r(x,D)f=f \in C^{\infty}(\Omega).
\end{equation}
Writing $\Delta_{\bsL}q(x,D)f=\diver_{\bsL^{*}}u$ for $u_{j}=L_{j}q(x,D)f$ we then get: 
\begin{equation}\nonumber
\diver_{\bsL^{*}}u-f=r(x,D)f
\end{equation}
for every  $f \in C^{\infty}(\Omega)$.

As application from the previous identity we present the following a priori estimates
{
\begin{proposition}\label{prop.gradmod}
Assume that the system of vector fields $L_{1},...,L_{n}$, $n \geq 2$, is linearly independent  and elliptic. Then for every point $x_{0} \in \Omega$ and $0<\beta<1$, there exist an open neighborhood $U \subset \Omega$ and a constant $C=C(U)>0$ such that, for all $\varphi\in\calD(U)$, one has:
\begin{equation}\label{jh}
\| \varphi \|_{{1-\beta,1}}:= \|J_{\beta-1}\varphi\|_{{1}}\leq C\|\nabla_{\bsL}\varphi\|_{{1}}.
\end{equation}
\end{proposition}
In the above statement, the operator $J_{\alpha}:=J_{\alpha}(x,D)$ for $\alpha>0$ is the pseudodifferential operator, called Bessel potential,  defined by 
$$J_{\alpha}f(x)=\int_{\R^{N}}e^{2\pi i x \cdot \xi}b(x,\xi)\hat{f}(\xi)d\xi, \quad f \in S'(\mathbb{R}^{N}),$$
where the symbol $b(x,\xi)=\langle \xi \rangle^{\alpha}:=(1+4\pi^{2}|\xi|^{2})^{-\alpha/2}$, independent of $x$, belongs to the H\"ormander class $S^{-\alpha}_{1,0}(\R^{N})$.
The operator $J_{-\alpha}$, usually denoted by $(1-\Delta)^{\alpha/2}$, allows us to introduce a nonhomogeneous fractional Sobolev space $W^{\alpha,p}(\R^{N})$ for $1\leq p <\infty$, defined
as the set of tempered distributions $u \in S'(\R^{N})$ such that $J_{-\alpha}u \in L^{p}(\R^{N})$, endowed with the norm $\|u\|_{{\alpha,p}}:=\|J_{-\alpha}u\|_{L^{p}}$.   As a consequence of the continuity property of the action of the Bessel potential on Lebesgue spaces (see for instance \cite[Theorem~2.5]{AH}), the inclusion $W^{\alpha,p}(\R^{N}) \subset L^{p}(\R^{N})$ is continuous for all $1\leq p<\infty$.
\begin{proof}
Let $h=\nabla_{\bsL}\varphi$. Thanks to identity \eqref{pseudo} we have 
$$J_{\beta-1}\varphi=p(x,D)h+r'(x,D)\varphi,$$
where $r'(x,D)=J_{\beta-1}r(x,D)$ is a regularizing operator and $p(x,D)=J_{\beta-1}(x,D)q_{1}(x,D)\diver_{\bsL^{\ast}}$ is a vector-valued pseudodifferential operator of negative order $-\beta$. As a consequence of Theorem \ref{l1} we have $\|p(x,D)h\|_{L^{1}}\lesssim \|\nabla_{\bsL}\varphi\|_{L^{1}}$, which implies:
 \begin{equation}\nonumber
\| J_{\beta-1}\varphi\|_{{1}}\leq C\|\nabla_{\bsL}\varphi\|_{L^{1}} +\|r'(x,D)\varphi\|_{{1}}.
\end{equation}
As the second term on the right side may be absorbed (see \cite[p.~798]{HP1}), shrinking the neighborhood if necessary, we obtain the estimate \eqref{jh}.
\end{proof}}
{
%\begin{remark}
The boundedness in $\LL^1$ norm of the pseudodifferential operators with negative order {follow from the}  integrability property of the kernel due itself to {a} pointwise control obtained in \cite{AH}.  %  We present the $L^{1}$ boundedness for pseudodifferential operators with negative order in the Theorem \ref{l1} .    
%\end{remark}
Another fundamental tool from pseudodifferential operators theory, inspired in the recent results obtained in \cite{HKP}, asserts that the embedding  $W_{c}^{1-\beta, 1}(B):=W^{1-\beta,1}(\R^{N})\cap \mathcal{E}'(B)$, where $B$ is a generic ball, {into} $L^{1}(\R^{N})$ is compact. These results are stated in the Appendix and will be proved there for sake of completeness.}

\section{Functions of bounded $\bsL$-variation}\label{sec1}

Throughout this section, we consider $L_1,\dots,L_n$ a system of complex vector fields with smooth coefficients on $\Omega$.

\subsection{Basic definitions; approximation and compactness}
Let $L^{1}_{c}(\Omega)$ be the linear space of all complex functions in $L^{1}(\Omega)$ whose support is a compact subset of $\Omega$.

The following definition of $\bsL$-variation of $g\in \LL^1_c(\Omega)$ recalls the classical definition of variation in case $n=N$ and {$L_j=\partial_{x_{j}}$ for each $j=1,2,\dots, N$.} It has been formulated for (real) vector fields by N. Garofalo and D. Nhieu \cite{GAROFALONHIEU}.
\begin{Definition}
Given $g\in\LL^1_c(\R^n)$ and $U\subseteq\Omega$ an open set, one calls the extended real number:
$$
\|D_\bsL g\|(U):=\sup\left\{\left|\int_\Omega g\,\overline{\diver_{\bsL^*} v}\right|: v\in C^\infty_c(\Omega,\bbC^n), \supp v\subseteq U, \|v\|_\infty\leq 1\right\},
$$
the \emph{(total) $\bsL$-variation of $g$ in $U$} and we let $\|D_\bsL g\|:=\|D_\bsL g\|(\Omega)$ in case there is no ambiguity on the open set $\Omega$.
We denote by $BV_{\bsL,c}(\Omega)$ the set of all $g\in \LL^1_c(\Omega)$ with $\|D_\bsL g\|<+\infty$.

Given $g\in BV_{\bsL,c}(\Omega)$, we denote by $D_\bsL g$ the unique $\bbC^n$-valued Radon measure satisfying:
\begin{equation}\label{eq.def-mes}
\int_\Omega g\,\overline{\diver_{\bsL^*} v}=\int_\Omega \bar{v}\cdot d [D_\bsL g],
\end{equation}
for all $v\in C^\infty_c(\Omega,\bbC^n)$. It is clear by definition that $\|D_\bsL g\|$ is also the total variation in $\Omega$ of $D_\bsL g$.
\end{Definition}
The next proposition allows us to define a vector-valued Radon measure $D_\bsL g$ for any $g\in BV_{\bsL,c}(\Omega)$.
\begin{Remark}Given $g\in BV_{\bsL,c}(\Omega)$, one has {$\supp D_\bsL g\subseteq \supp g$}. Indeed, given $x\in \Omega\setminus\supp g$, find a radius $r>0$ for which one has $B(x,r)\subseteq \Omega\setminus \supp g$. It is clear according to (\ref{eq.def-mes}) that for any $v\in C^\infty_c(B(x,r),\bbC^n)$ we then have $D_\bsL g(v)=0$. Hence we also get $D_\bsL g(v)=0$ for all $v\in C_c(B(x,r),\bbC^n)$, which ensures that one has $x\notin\supp D_\bsL g$ and finishes to show the inclusion $\supp D_\bsL g\subseteq\supp g$.
\end{Remark}
\begin{Remark}\label{rmk.sci}
It follows readily from the previous definition that, as in the classical case, if $(g_i)\subseteq BV_{\bsL,c}(\Omega)$ converges in $\LL^1$ to $g\in \LL^1_c(\Omega)$, one then has $g\in BV_{\bsL,c}(\Omega)$ and:
$$
\|D_\bsL g\|\leq\liminfe_{i}\|D_\bsL g_i\|.
$$
We shall refer to this in the sequel as the \emph{lower semi-continuity} of the $\bsL$-variation.
\end{Remark}

We say that a sequence $(f_{i})_{i}$ of functions with complex values defined on open set $\Omega \subset \R^{N}$ is \textit{compactly supported in $\Omega$} if there is a compact set $K \subset\subset \Omega$ such that one has $\supp f_i \subseteq K$ for every $i$.

We shall make an extensive use of the following concept of convergence.
\begin{Definition}
Given $g \in \LL^1_{c}(\Omega)$ and a sequence $(\varphi_i)_{i}\subseteq\calD(\Omega)$ we shall write $\varphi_i \twoheadrightarrow g $ in case the following conditions hold:
\begin{enumerate}
\item[(i)] $(\varphi_{i})$ converges to $g$ in $\LL^{1}$ norm;
\item[(ii)] $(\varphi_i)$ is compactly supported in $\Omega$;
\item[(iii)] $\sup_i\| \nabla_\bsL \varphi_{i}\|_1<+\infty$.
\end{enumerate}
\end{Definition}

Using a Friedrich's type decomposition due to N. Garofalo and D. Nhieu \cite[Lemma~A.3]{GAROFALONHIEU} in the real case, we obtain an analogous result, in $BV_{\bsL,c}$, to the standard approximation theorem for $BV_c$ functions.
\begin{lemma}\label{lemmaaprox}\label{lem.approx}
Assume that $L_1,\dots, L_n$ have locally Lipschitz coefficients.
For any $g \in BV_{\bsL,c}(U)$, there exists a sequence $\left\{ \varphi_{i} \right\}_{i} \subset \calD(U) $ such that one has $\varphi_i\twoheadrightarrow g$ and, moreover:
$$
\|D_\bsL g\|=\lim_i\| \nabla_\bsL \varphi_{i}\|_1.
$$
\end{lemma}
\begin{proof}
Fix $\eta\in\calD(\R^n)$ a radial function with nonnegative values, satisfying $\supp \eta\subseteq B[0,1]$ and $\int_{\R^n} \eta=1$, and, for each $\epsilon>0$, define $\eta_\epsilon\in\calD(\R^n)$ by $\eta_\epsilon(x):=\epsilon^{-N}\eta(x/\epsilon)$.

Fix now $g\in BV_{\bsL,c}(\Omega)$ and define for $0<\epsilon<\dist(\supp g, \complement\Omega)$ a function $g_\epsilon\in\calD(\Omega)$ by the formula:
$$
g_\epsilon:=\eta_\epsilon* g.
$$
For each $i=1,\dots, n$, denote by $D_{L_i}g$ the compactly supported distribution defined by:
$$
D_{L_i}g (\varphi):=\int_\Omega g\, \overline{L_i^*\varphi},
$$
let $D_\bsL g$ denote the vector-valued distribution $(D_{L_1}g, \dots, D_{L_n}g)$
and observe that according to N. Garofalo and D. Nhieu \cite[Lemma~A.3]{GAROFALONHIEU}, one can write:
\begin{equation}\label{eq.friedrichs}
\nabla_\bsL(\eta_\epsilon*g)=\eta_\epsilon *(D_{\bsL} g)+H_{\epsilon}(g),
\end{equation}
where also $\|H_{\epsilon} (g)\|_1\to 0$, $\epsilon\to 0$.

Fix now $v\in C^\infty_c(\Omega, \bbC^n)$ a smooth vector field satisfying $\|v\|_\infty\leq 1$ and compute:
\begin{multline*}
\left|\int_\Omega \nabla_{\bsL} (\eta_\epsilon*g)\cdot \bar{v}\right|\leq \left|\sum_{j=1}^n D_{L_i} g (\eta_\epsilon*v_i)\right|+\|v\|_\infty \|H_\epsilon(g)\|_1\\
=\left|\int_\Omega g\,\overline{\diver_{\bsL^*} (\eta_\epsilon*v)}\right|+\|H_\epsilon(g)\|_1\leq \|D_\bsL g\|+\|H_\epsilon(g)\|_1.
\end{multline*}
We hence get, by duality:
$$
\|\nabla_{\bsL} (\eta_\epsilon* g)\|_1\leq\|D_\bsL g\|+\|H_\epsilon(g)\|_1,
$$
and the result follows from the aforementioned property of $H_\epsilon(g)$ when $\epsilon$ approaches $0$.
\end{proof}

The following proposition is a compactness result in $BV_{\bsL}$.
\begin{Proposition}\label{prop.compacite}
Assume that the open set $U\subseteq \Omega$ supports a Sobolev-Gagliardo-Nirenberg inequality of type (\ref{SGN}) as well as an inequality of type (\ref{jh}) for some $\epsilon>0$. If $(g_i)\subseteq BV_{\bsL,c}(U)$ is compactly supported in $U$ and if moreover one has:
$$
\sup_{i} \|D_\bsL g_i\|<+\infty,
$$
then there exists $g\in BV_{\bsL,c}(U)$ and a subsequence $(g_{i_k})\subseteq (g_i)$ converging to $g$ in $\LL^1$ norm.
\end{Proposition}
\begin{proof}
Choose a compact set $K\subset\subset U$ for which one has $\supp g_i\subseteq K$ for all $i$, and let $\chi\in\calD(U)$ be such that $\chi_K\leq\chi\leq 1$ on $U$. Choose also, according to Lemma~\ref{lem.approx}, a sequence $(\varphi_i)\subseteq\calD(U)$ satisfying the following conditions for all $i$:
$$
\|g_i-\varphi_i\|_1\leq 2^{-k}\quad\text{and}\quad \|\nabla_\bsL\varphi_i\|_1\leq \|D_\bsL g_i\|+1.
$$
Define now, for each $i$, $\psi_i:=\varphi_i\chi$ and compute using H\"older's inequality together with (\ref{SGN}):
$$
\|\nabla_\bsL \psi_i\|_1\leq \|\varphi_i\|_{N/N-1} \|\nabla_\bsL \chi\|_N + \|\nabla_\bsL \varphi_i\|_1\leq (C\|\nabla\chi\|_N +1) \|\nabla_\bsL\varphi_i\|_1.
$$
We hence have $\sup_i \|\nabla_\bsL \psi_i\|_1<+\infty$ while it is clear that $(\psi_i)$ is compactly supported and satisfies $\|g_i-\psi_i\|_1\to 0$, $i\to\infty$.

Now fix $0<\beta<1$ and observe that the sequence $(\psi_i)_{i}$ also satisfies, according to (\ref{jh}):
$$
\sup_i \|\psi_i\|_{1-\beta,1}=\sup_i \|J_{\beta-1} \psi_i\|_1\leq C\sup_i \|\nabla_\bsL \psi_i\|_1<+\infty.
$$
It hence follows from the compactness of the inclusion of $W^{1-\beta,1}_c(U)\subset \subset\LL^1(U)$ (see Theorem~\ref{thm.compact-beta} in Appendix) that there exists $g\in \LL^1(U)$ and a subsequence $(\psi_{i_k})\subseteq (\varphi_i)$ converging to $g$ in $\LL^1(U)$. On the other hand it is clear that one has $\supp g\subseteq K$ as well as $\psi_{i_k}\twoheadrightarrow g$. We hence have, by lower semicontinuity:
$$
\|D_\bsL g\|\leq\liminfe_{k} \|\nabla_\bsL \psi_{i_k}\|_1,
$$
which ensures that one has $g\in BV_{\bsL,c}(U)$.

\end{proof}
\begin{Remark}
According to Theorem~\ref{thmHP} and Proposition~\ref{prop.gradmod}, we see that if one assumes $L_1,\dots, L_n$ to be everywhere linearly independent and elliptic, each point $x_0\in\Omega$ is contained a neighborhood $U\subseteq\Omega$ satisfying the hypotheses of the previous proposition. \end{Remark}

\subsection{A Sobolev-Gagliardo-Nirenberg inequality in $BV_\bsL$} As announced we get the following result:

\begin{proposition}
Assume that the system of vector fields $L_{1},...,L_{n}$, $n \geq 2$, is linearly independent  and elliptic. Then every point $x_{0} \in \Omega$ is contained in an open neighborhood $U \subset \Omega$ such that the inequality:
\begin{equation}\label{dl}
\|g\|_{{N/N-1}}\leq C\|D_{\bsL}g\|,
\end{equation}
holds for all $g\in BV_{\bsL,c}(U)$, where $C=C(U)>0$ is a constant depending only on $U$.
\end{proposition}
\begin{proof}
Fix $x_{0} \in \Omega$. It follows from Theorem \ref{thmHP} that there exists a neighborhood $U \subset \Omega$ of $x_{0}$ and $C=C(U)>0$ such that, for all $\varphi\in\calD(U)$, one has:
\begin{equation}\nonumber
\|\varphi\|_{{N/N-1}}\leq C \| \nabla_{\bsL} \varphi\|_{{1}}.
\end{equation}
Then given $g \in BV_{\bsL,c}(U)$ consider the sequence $\left\{ \varphi_{i} \right\} \subset \calD(U)$ satisfying (i)-(iii) by Lemma \ref{lemmaaprox}. As a consequence of Fatou Lemma and the previous estimate we conclude that
\begin{equation}\nonumber
\|g\|_{{N/N-1}}\leq \limsupe_{i \rightarrow \infty}\|\varphi_{i}\|_{{N/N-1}} 
\leq C  \limsupe_{i \rightarrow \infty}\|\nabla_{\bsL}\varphi_{i}\|_{{1}} 
\leq C' \|D_{\bsL}g\|.
\end{equation}
 The proof is complete.
\end{proof}

\begin{Remark}
The converse of proposition is true, namely if the inequality \eqref{dl} holds then the system must be elliptic on $U$ (see \cite{HP1} for details).
\end{Remark}

\section{$\bsL$-charges and their extensions to $BV_{\bsL,c}$}\label{sec2}

%\begin{definition}
%The distribution $F:\Omega \rightarrow \Lambda^{n}\R^{n}$ is a $\bsL$-strong charge if for every $\epsilon>0$ and $R>0$, there exists $C>0$ such that if $\varphi \in C_{c}^{\infty}(B_{R})$,
%\begin{equation}\label{lstrong}
%\left| \int_{\Omega} F \wedge  \varphi \right| \leq C\|\varphi\|_{L^{1}}+\epsilon\|\nabla_{\bsL} \varphi\|_{L^{1}}.
%\end{equation}
%\end{definition}

We now get back to the original problem of finding, locally, a continuous solution to (\ref{diverl}). 
\subsection{$\bsL$-fluxes and $\bsL$-charges} Distributions which allow, in an open set $\Omega$, to solve continuously (\ref{diverl}), will be called \emph{$\bsL$-fluxes}.

\begin{definition}
A distribution $F \in \calD'{(\Omega)}$ is called an \emph{$\bsL$-flux in $\Omega$} if the equation \eqref{diverl} % $\diver_{\bsL^{\ast}}v=F$
has a continuous solution in $\Omega$, \emph{i.e.} if there exists $v \in C(\Omega,\mathbb{C}^{n})$ such that one has, for all $\varphi\in\calD(\Omega)$:
\begin{equation}\label{eq.flux}
F(\varphi)=\int_{\Omega}\bar{v}\cdot {\nabla_{\bsL}\varphi}, \quad \forall\; \varphi \in \calD(\Omega).
\end{equation}
\end{definition}

$\bsL$-fluxes satisfy the following continuity condition.
\begin{lemma}
If $F$ is an $\bsL$-flux then $\lim_{i} F(\varphi_{i})=0$ for every sequence $( \varphi_{i})_{i}\subseteq \calD(\Omega)$ verifying $\varphi_i\twoheadrightarrow 0$.
\end{lemma}
\begin{proof}
Let $F$ be an $\bsL$-flux and let $v\in C(\Omega,\bbC^n)$ be such that (\ref{eq.flux}) holds.
Fix a sequence $(\varphi_i)_{i}\subseteq\calD(\Omega)$ verifying $\varphi_i\twoheadrightarrow 0$, let $c:=\sup_i\|\nabla_\bsL\varphi_i\|_1<+\infty$ and choose a compact set $K\subset\subset\Omega$ for which one has $\supp\varphi_i\subseteq K$ for all $i$.

Fix now $\epsilon>0$. According to Weierstrass' approximation theorem, choose a vector field $w\in C^\infty_c(\Omega,\bbC^n)$ for which one has $\sup_K |v-w|\leq\epsilon$ and compute, for all $i$:
$$
\left| F(\varphi_i)\right|\leq \left|\int_\Omega (\bar{v}-\bar{w})\cdot {\nabla_\bsL \varphi_i}\right|+
\left|\int_\Omega \bar{w}\cdot {\nabla_\bsL \varphi_i}\right|\leq\epsilon \|\nabla_\bsL\varphi_i\|_1+\left| \int_\Omega {\varphi_i} \,\overline{\diver_{\bsL^*}w}\right|\leq c\epsilon+ \|\diver_{\bsL^*} w\|_\infty \|\varphi_i\|_1.
$$
We hence get $\limsupe_i |F(\varphi_i)|\leq c\epsilon$,
and the result follows for $\epsilon>0$ is arbitrary.
\end{proof}

The property above suggest the following definition of linear functionals associated to $\bsL$. 
\begin{definition}
A linear functional $F: \calD(\Omega) \rightarrow \mathbb{C}$ is called an $\bsL$-charge in $\Omega$  if $\lim_{i} F(\varphi_{i})=0$ for every sequence $( \varphi_{i})_{i}\subseteq \calD(\Omega)$ satisfying $\varphi_i\twoheadrightarrow 0$. The linear space of all $\bsL$-charges in $\Omega$ is denoted by $CH_{\bsL}(\Omega)$.
\end{definition}

The following characterization of $\bsL$-charges will be useful in the sequel.
\begin{proposition}
If $F:\calD(\Omega) \rightarrow \mathbb{C}$ is a linear functional, then the following properties are equivalent
\begin{enumerate}
\item[(i)] $F$ is an $\bsL$-charge,
\item[(ii)] for every $\epsilon>0$ and each compact set $K \subset \subset \Omega$ there exists $\theta>0$ such that, for any $\varphi\in\calD_K(\Omega)$, one has:
\begin{equation}\label{eq.ch}
|F(\varphi)|\leq \theta\|\varphi\|_{L^{1}}+\epsilon\|\nabla_{\bsL}\varphi\|_1.
\end{equation}
\end{enumerate}
\end{proposition}
\begin{proof}
We proceed as in \cite[Proposition~2.6]{PP}.\

Since (ii) implies trivially (i), it suffices to show that the converse implication holds. To that purpose, assume (i) holds, \emph{i.e.} suppose that $F$ is an $\bsL$-charge. Fix $\epsilon>0$ and a compact set $K \subset \subset \Omega$. By hypothesis, there exists $\eta>0$ such that for every $\varphi \in \calD_K(\Omega)$ satisfying $\| \varphi \|_{{1}}\leq \eta$ and $\|D_{\bsL}\varphi\|_{{1}}\leq 1$, we have $|F(\varphi)|\leq \epsilon$. We now define $\theta :=\epsilon / \eta$.

Fix now $\varphi\in\calD_K(\Omega)$ and assume by homogeneity that one has $\|\nabla_{\bsL}\varphi\|_1=1$. If moreover one has $\|\varphi\|_{{1}}\leq \eta$, then one computes $|F(\varphi)|\leq \epsilon=\epsilon \|\nabla_{\bsL}\varphi\|_1$. If on the contrary we have $\|g\|_{L^{1}}> \eta$, we define $\tilde{\varphi}=\varphi \eta/\|\varphi\|_{{1}}$. We then have $\|\tilde{\varphi}\|_{{1}}= \eta$ as well as $\|\nabla_{\bsL}\tilde{\varphi}\|_{{1}}< 1$, and hence also $|F(\tilde{\varphi})|\leq \epsilon$; this yields finally $|F(\varphi)|=\|\varphi\|_{{1}}|f(\tilde{\varphi})|/\eta\leq \epsilon \|\varphi\|_{{1}}/\eta=\theta \|\varphi\|_{{1}}$.
\end{proof}

As we shall see now, $\bsL$-charges can be extended in a unique way to linear forms on $BV_{\bsL,c}$.
\begin{Proposition}
An $\bsL$-charge $F$ in $\Omega$ extends in a unique way to a linear functional $\tilde{F}:BV_{\bsL,c}(\Omega)\to\bbC$ satisfying the following property: for any $\epsilon>0$ and each compact set $K\subset\subset\Omega$, there exists $\theta>0$ such that for any $g\in BV_{\bsL,K}(\Omega)$ one has:
\begin{equation}\label{eq.ch2}
|\tilde{F}(g)|\leq\theta\|g\|_1+\epsilon\|D_\bsL g\|.
\end{equation}
\end{Proposition}
\begin{proof}
Given $g\in BV_{\bsL,c}(\Omega)$, fix $(\varphi_i)_{i}\subseteq\calD(\Omega)$ satisfying $\varphi_i\twoheadrightarrow g$ and observe that it follows from (\ref{eq.ch}) that $(F(\varphi_i))_i$ is a Cauchy sequence of complex numbers whose limit does not depend on the choice of sequence $(\varphi_i)\subseteq\calD(\Omega)$ satisfying $\varphi_i\twoheadrightarrow g$. We hence define $\tilde{F}(g):=\lim_i F(\varphi_i)$. It now follows readily from (\ref{eq.ch}) and Lemma~\ref{lem.approx} that $\tilde{F}$ satisfies the desired property.
\end{proof}
\begin{Remark}\label{rmk.conv-ch}
If $\tilde{F}:BV_{\bsL,c}(\Omega)\to\bbC$ extends the $\bsL$-charge $F$, it is easy to see from the previous proposition that for any compactly supported sequence $(g_i)_{i}\subseteq BV_{\bsL,c}(\Omega)$ satisfying $g_i\to 0$, $i\to\infty$ in $L^1(\Omega)$ and $\sup_i \|D_\bsL g_i\|<+\infty$, one has $F(g_i)\to 0$, $i\to\infty$.
\end{Remark}
From now on, we shall identify any $\bsL$-charge with its extension to $BV_{\bsL,c}$ and use the same notation for the two linear forms.

\subsection{Two important examples of $\bsL$-charges} Let us define two important classes of $\bsL$-charges.
\begin{Example}\label{ex.flux}
In case $F$ is the $\bsL$-flux associated to $v\in C(\Omega, \bbC^n)$ according to (\ref{eq.flux}), its unique extension to $BV_{\bsL,c}(\Omega)$ is the $\bsL$-charge:
$$
\Gamma(v):BV_{\bsL, c}(\Omega)\to\bbC, g\mapsto \int_{\Omega} \bar{v}\cdot d\left[{D_\bsL g}\right].
$$

To see this, fix $g\in BV_{\bsL,c}(\Omega)$ together with a sequence $(\varphi_i)_{i}\subseteq\calD(\Omega)$ satisfying (i)-(iii) in Lemma~\ref{lem.approx} and choose $\supp g\subseteq K\subset\subset \Omega$ a compact set for which one has $\supp _i\subseteq K$ for all $i$.
Given $\epsilon>0$, choose $w\in C^\infty_c(\Omega,\bbC^n)$ a smooth vector field satisfying $\sup_K |v-w|\leq\epsilon$ and compute:
$$
\left|\Gamma(v)(g)-\int_\Omega \bar{v}\cdot d\left[{D_\bsL g}\right]\right|=\lim_i \left| \int_\Omega \bar{v}\cdot {\nabla_\bsL \varphi_i}-\int_\Omega \bar{v}\cdot d\left[{D_\bsL g}\right]\right|.
$$
On the other hand we have for all $i$:
\begin{multline*}
\left| \int_\Omega \bar{v}\cdot {\nabla_\bsL \varphi_i}-\int_\Omega \bar{v}\cdot d\left[ {D_\bsL g}\right]\right|\leq \left| \int_\Omega (\bar{v}-\bar{w})\cdot {\nabla_\bsL \varphi_i}\right|+\left| \int_\Omega (\bar{v}-\bar{w})\cdot d\left[ {D_\bsL g}\right]\right|\\+\left| \int_\Omega \bar{w}\cdot {\nabla_\bsL \varphi_i}-\int_\Omega \bar{w}\cdot d\left[ {D_\bsL g}\right]\right|
\leq \epsilon \|\nabla_\bsL\varphi_i\|_1+\epsilon \|D_\bsL g\| +\left| \int_\Omega {\varphi_i}\, \overline{\diver_{\bsL^*} w}-\int_\Omega \bar{w}\cdot d\left[ {D_\bsL g}\right]\right|.
\end{multline*}
Using the properties of $(\varphi_i)_{i}$ and Lebesgue's dominated convergence theorem, we thus get:
$$
\lim_i \left| \int_\Omega \bar{v}\cdot {\nabla_\bsL \varphi_i}-\int_\Omega \bar{v}\cdot d\left[ {D_\bsL g}\right]\right|\leq 2\epsilon \|D_\bsL g\|+\left| \int_\Omega {g}\,\overline{ \diver_{\bsL^*} w}-\int_\Omega \bar{w}\cdot d\left[ {D_\bsL g}\right]\right|=2\epsilon \|D_\bsL g\|,
$$
according to (\ref{eq.def-mes}). The result follows, for $\epsilon>0$ is arbitrary.
\end{Example}
\begin{Example}\label{ex.lambda}
Assume that $U$ supports a Sobolev-Gagliardo-Nirenberg inequality of type (\ref{dl}) for $BV_{\bsL}$ functions in $U$. Define then, for any $f\in L^N(U)$, a map $\Lambda(f):BV_{\bsL,c}(U)\to\bbC$ by:
$$
\Lambda(f)(g):=\int_U \bar{f}{g}.
$$
Fix $\epsilon>0$ and choose $\theta>0$ large enough for $\int_{\{|f|>\theta\}} |f|^N\leq \epsilon^{N}$ to hold. We then compute:
\begin{eqnarray*}
{\int_{\Omega}|\bar{f}{g}|}&\leq &{\theta \int_{|f|\leq \theta}|g|+\int_{|f|> \theta} |fg|},\\
&\leq &{\theta \|g\|_{{1}}+\left(\int_{\{|f|>\theta\}} |f|^N\right)^{\frac 1N}\|g\|_{{N/N-1}}},\\
&\leq &{\theta \|g\|_{L^{1}}+\epsilon\|D_{\bsL}g\|}.
\end{eqnarray*}
for appropriated choice of $\theta$. Hence $\Lambda(f)$ defines an $\bsL$-charge.
\end{Example}
\begin{Remark}\label{rmk.ODE}
It is easy to see that for any $x_0=(x_0^1,\dots, x_0^N)\in\Omega$, there exists an open set $x_0\in U\subseteq\Omega$ such that one has $\Lambda[\calD(U)]\subseteq\Gamma[C^\infty(U,\bbC^n)]$.

%\textcolor{red}{Explain here how one can solve locally $\Delta_\bsL u=f$ in $C^\infty$.}

{Given $\varphi\in\calD(U)$, thanks to the local solvability of the elliptic equation\eqref{deltal} (see \cite[Corollary 4.8]{GS}), there exists $u\in C^\infty(U)$ a smooth solution to $\Delta_\bsL u=\varphi$ in $U$.  Let $v:=\nabla_\bsL u$.}
%Now given $\varphi\in\calD(U)$, write $u\in C^\infty(U)$ a smooth solution to $\Delta_\bsL u=\varphi$ in $U$, \textcolor{red}{thanks to $\Delta_{\bsL}$ is elliptic and then has local solvability (see Cor 4.8 in \cite{GS})}. Let $v:=\nabla_\bsL u$.
This yields, for any $g\in BV_{\bsL,c}(U)$:
$$
\Lambda(\varphi)(g)=\int_U \bar{\varphi} {g}=\int_U g \,\overline{\diver_{\bsL^*} v}=\int_U \bar{v}\cdot d\left[{D_\bsL g}\right]=\Gamma(v)(g),
$$
for we could, in the computation above, replace $v$ by $v\chi$ where $\chi\in\calD(U)$ satisfies $\chi=1$ in a neighborhood of $\supp g$.
\end{Remark}

It turns out that a linear functional on $BV_{\bsL,c}$ is an $\bsL$-charge if and only if it is continuous with respect to some locally convex topology on $BV_{\bsL,c}$.

\subsection{Another characterization of $\bsL$-charges}
In the sequel, a \emph{locally convex space} means a Hausdorff locally convex topological vector space. For any family  $\mathcal{A}$ of sets and any set $E$ we denote $\mathcal{A}  \lefthalfcup  E :=\left\{ A \cap E: A \in \mathcal{A} \right\}$. Following \cite[Theorem~3.3]{PMP} we define the following topology on $BV_{\bsL,c}(\Omega)$ (note that this result remains valid in the complex framework).

\begin{definition}
Let $\bsT_{\bsL}$ be the unique locally convex topology on $BV_{\bsL,c}(\Omega)$ such that
\begin{enumerate}
\item[(a)] $\bsT_{\bsL} \lefthalfcup BV_{\bsL,K,\lambda} \subseteq \bsT_{\LL^{1}} \lefthalfcup BV_{\bsL,K,\lambda}$ for all $K \subset \subset \Omega$ and $\lambda>0$ where we let:
$$BV_{\bsL,K,\lambda}=\left\{ g \in BV_{\bsL,c}(\Omega):\;\text{supp}\;g\subseteq K, \|D_{\bsL}g\|\leq \lambda \right\},$$
and where $\bsT_{\LL^1}$ is the $\LL^1$-topology;
\item[(b)] for every locally convex space $Y$, a linear map $f:(BV_{\bsL,c};\bsT_{\bsL}) \rightarrow Y$ is continuous if only if $ f\upharpoonright{BV_{\bsL,K,\lambda}}$ is $\LL^{1}$ continuous for all $K \subset \subset \Omega$ and $\lambda>0$ .
\end{enumerate}
\end{definition}

$\bsL$-charges are the $\bsT_\bsL$-continuous linear functionals, as it readily follows from Remark~\ref{rmk.conv-ch}.
\begin{Proposition}
A linear functional $F:BV_{\bsL,c}(\Omega)\to\bbC$ is an $\bsL$-charge if and only if it is $\bsC_\bsL$-continuous.
\end{Proposition}

We now turn to proving the key result for obtaining Theorem~\ref{mainthm}.
\section{Towards Theorem~\ref{mainthm}}\label{sec3}

Throughout this section, we assume that $L_1,\dots, L_n$ is a system of linearly independent vector fields in $\Omega$, and that the open set $U\subseteq\Omega$ supports inequalities of type (\ref{SGN}) and (\ref{jh}); we also assume that one has $\Lambda[\calD(U)]\subseteq \Gamma[C(U,\bbC^n)]$.
\begin{Remark}
It follows from Theorem~\ref{thmHP}, Proposition~\ref{prop.gradmod} and Remark~\ref{rmk.ODE} that for any $x_0\in \Omega$, one can find an open neigborhood $U$ of $x_0$ in $\Omega$ satisfying all the above assumptions.
\end{Remark}

Our intention is to prove the following result.

\begin{Theorem}
If $F:BV_{\bsL,c}(U)\to\bbC^n$ is an $\bsL$-charge in $U$, then there exists $v\in C(U,\bbC^n)$ for which one has $F=\Gamma(v)$, \emph{i.e.} such that one has, for any $g\in BV_{\bsL,c}(U)$:
$$
F(g)=\int_U \bar{v}\cdot d\left[{D_\bsL g}\right].
$$
\end{Theorem}
To prove this theorem, we have to show that the map
$$
\Gamma: C( U ,\mathbb{C}^{n}) \longrightarrow CH_{\bsL}( U ), v\mapsto \Gamma(v),
$$
is surjective. In order to do this, we endow $C(U,\bbC^n)$ with the usual Fr\'echet topology of uniform convergence on compact sets, and $CH_\bsL ( U )$ with the Fr\'echet topology associated to the family of seminorms $(\|\cdot\|_K)_K$ defined by:
$$\|F\|_{K} := \sup\left\{ | F(g) |: g \in BV_{\bsL,K}( U ), \|D_\bsL g\|\leq 1\right\},$$
where $K$ ranges over all compact sets $K\subset\subset U $. The surjectivity of $\Gamma$ will be proven in case we show that $\Gamma$ is continuous and verifies the following two facts:
\begin{enumerate}
\item[(a)] $\Gamma[C( U ,\mathbb{C}^{n})]$ is dense in $CH_{\bsL}( U )$.
\item[(b)] $\Gamma^{\ast}[CH_{\bsL}( U )^{\ast}]$ is sequentially closed in the strong topology of $C( U ,\bbC^n)^*$.
\end{enumerate}
Indeed, it will then follow from the Closed Range Theorem \cite[Theorem~8.6.13]{EDW} together with\cite[Proposition~6.8]{PMP} and (b) that $\Gamma[C( U ,\mathbb{C}^{n})]$ is closed in $CH_{\bsL}( U )$. Using (a) we shall then conclude that one has:$$\Gamma[C( U ,\mathbb{C}^{n})]=CH_{\bsL}( U ),$$
\emph{i.e.} that $\Gamma$ is surjective.

The strategy of the proof of (b) follow the lines of De Pauw and Pfeffer's proof in \cite{PP}. For the proof of (a), however, the proof presented below is slightly different from their approach; we namely manage to avoid an explicit smoothing process and choose instead to use an abstract approach similar to the one used in \cite{M} in order to solve the equation $d\omega=F$.

Let us start by showing that $\Gamma$ is continuous.
\begin{lemma}
The map $\Gamma: C( U ,\mathbb{C}^{n}) \longrightarrow CH_{\bsL}( U )$ is linear and continuous.
\end{lemma}
\begin{proof}
Indeed given a compact set $K\subset\subset  U $ and $g \in BV_{K,\bsL}( U )$ we have:
$$ |\Gamma(v)(g) | = \left|  \int_{ U } \bar{v}\cdot d\left[{D_{\bsL}g}\right] \right| \leq \| D_{\bsL}g \|\|v\|_{\infty,K},   $$
which implies  $\|\Gamma(v)\|_{K}\leq \|v\|_{\infty,K}$.
\end{proof}

First we have to identify the dual space $CH_{\bsL}( U )^{*}$.
\subsection{Identifying the dual space $CH_\bsL( U )^*$}
The following result is the identification we need.
\begin{proposition}
The map $\Phi: BV_{\bsL,c}( U ) \longrightarrow CH_{\bsL}( U )^{*}$ given by $\Phi(g)(F):=F(g)$ is a linear bijection.
\end{proposition}

The proof of the previous proposition is quite delicate. We shall proceed in several steps which will be interesting as such.\\

First let us check that $\Phi$ is well defined. In fact, given $K \subset \subset  U $ and $g\in BV_{\bsL,K}(U)$ we have
$$ \left|\Phi(g)(F) \right|= \left|  F(g) \right|\leq \|D_{\bsL}g\|\|F\|_{K},$$
according to the definition of $\|\cdot\|_{K}$. Hence $\Phi(g)$ is continuous and $\Phi(g) \in CH_{\bsL}( U )^{\ast}$. \\

To show that $\Phi$ is injective, let $g \in BV_{\bsL,c}( U )$ be such that $\Phi(g)=0$. Then for any $B \subset  U $ measurable and bounded we have:
\begin{equation}\nonumber
\int_{B}{g}=\int_{ U }\chi_{B}{g}=\Lambda(\chi_{B})(g)=\Phi(g)[\Lambda(\chi_{B})]=0.
\end{equation}
Thus $g=0$ a.e. in $ U $, which implies that $\Phi$ injective.\\

The next step is to prove that $\Phi$ is surjective. To show this property we shall define a right inverse for $\Phi$, called $\Psi$.\\

Let $\Psi: CH_{\bsL}( U )^{\ast} \longrightarrow \calD' ( U )$ be defined by:
\begin{equation}
\Psi(\alpha)[\varphi] :=\alpha[\Lambda(\varphi)].
\end{equation}
We claim that $\Psi$ is well defined, \emph{i.e.} that for $\alpha \in CH_{\bsL}( U )^{\ast}$, we have $\Psi(\alpha) \in BV_{\bsL,c}( U )$. Indeed, given $\alpha \in CH_{\bsL}( U )^{\ast}$ there exist $C>0$ and $K \subset \subset  U $ such that for all $F \in CH_{\bsL}( U )$ we have $|\alpha(F)|\leq C\|F\|_{K}$. In particular, for every $\varphi \in \calD( U )$ we have:
\begin{eqnarray*}
{|\Psi(\alpha)(\varphi)| }&\leq &C \|\Lambda(\varphi)\|_{K},\\
& \leq &C{\sup \left\{ |\Lambda(\varphi)(g)| :  g \in BV_{\bsL,K}( U ),  \|D_{\bsL}g\|\leq 1  \right\}},\\
&\leq &{ C \sup \left\{ \int_{ U }|\bar{\varphi}g|: g \in BV_{\bsL,K}( U ),  \|D_{\bsL}g\|\leq 1  \right\}},\\
&\leq &{C\|\varphi\|_{{N}} \sup \left\{ \|g\|_{{N/N-1}}  : g \in BV_{\bsL,K}( U ),  \|D_{\bsL}g\|\leq 1  \right\}},\\
& \leq &{C'\|\varphi\|_{{N}}  \sup \left\{ \|D_{\bsL}g\|:  g \in BV_{\bsL,K}( U ),  \|D_{\bsL}g\|\leq 1  \right\}},\\&\leq & C'\|\varphi\|_N,
\end{eqnarray*}
which implies that $\Psi(\alpha) \in L^{\frac{N}{N-1}}( U )$ by Riesz Representation theorem. Note that if $\varphi \in \calD( U )$  satisfies $(\supp\varphi )\cap K=\emptyset$ then one has
$|\Psi(\alpha)[\varphi]|=0$, which implies $\supp[\Psi(\alpha)] \subset K$.
Moreover, for any $v \in C^{\infty}_{c}( U ,\bbC^{n})$ we have:
\begin{eqnarray*}
 |\Psi(\alpha)[\diver_{\bsL^{\ast}}v]|&=& |\alpha[\Lambda(\diver_{\bsL^{\ast}}v)]|,\\ & \leq&{ C \|\Lambda(\diver_{\bsL^{\ast}}v)\|_{K}},\\
&\leq & C{\sup \left\{ \left|\int_{ U }g\,\overline{\diver_{\bsL^{\ast}}v}\right|: g \in BV_{\bsL,K}( U ),  \|D_{\bsL}g\|\leq 1  \right\}},\\
& \leq & C \sup \left\{ \|D_{\bsL}g\| \|v\|_{\infty}  : g \in BV_{\bsL,K}( U ),  \|D_{\bsL}g\|\leq 1  \right\},\\
&\leq  &C{\sup \|v\|_{\infty},  }
\end{eqnarray*}
so that one has $\Psi(\alpha) \in BV_{\bsL,c}( U )$. 

\begin{lemma}\label{lem.inv} The maps $\Phi$ and $\Psi$ defined above are inverses, \emph{i.e.} we have:
\begin{enumerate}
\item[(i)] $\Psi \circ \Phi=Id_{BV_{\bsL,c}( U )}$;
\item[(ii)] $\Phi \circ \Psi=Id_{CH_{\bsL}( U )^{*}}$ (in particular, $\Phi$ is surjective).
\end{enumerate}
\end{lemma}
In order to prove the previous lemma, we shall need some observations concerning the polar sets of some neighborhoods of the origin in $CH_{\bsL}( U )$. First, observe that the family of all sets $V(K,\epsilon)$ (where $K$ ranges over all compact subsets of $ U $, and $\epsilon$ over all positive real numbers) defined by:
$$
V(K,\epsilon):=\{F\in CH_\bsL( U ):\|F\|_K\leq\epsilon\},
$$
is a basis of neighborhoods of the origin in $CH_\bsL ( U )$.
\begin{Claim}\label{cl.w*}
Fix $K\subset\subset U $ a compact set and a real number $\epsilon>0$. For any $\alpha\in V(K,\epsilon)^\circ$, one has:
\begin{itemize}
\item[(i)] $\supp\Psi(\alpha)\subseteq K$;
\item[(ii)] $\|D_\bsL \Psi(\alpha)\|\leq\frac{1}{\epsilon}$.
\end{itemize}
\end{Claim}
\begin{proof}
To prove (i), assume that $\varphi\in\calD(\bbC^n)$ satisfies $K\cap \supp\varphi=\emptyset$. Then, we get for $\lambda>0$:
$$
\|\lambda\Lambda(\varphi)\|_K=\sup\left\{\lambda \left|\int_ U  \bar{\varphi} g\right|:g\in BV_{\bsL, K,\lambda}( U ), \|D_\bsL g\|\leq 1\right\}=0.
$$
In particular this yields $\lambda \Lambda(\varphi)\in V(K,\epsilon)$. We hence obtain:
$$
\lambda |\alpha[\Lambda(\varphi)]|=|\alpha[\lambda\Lambda(\varphi)]|\leq 1,
$$
for any $\lambda>0$. Since $\lambda>0$ is arbitrary, this implies that one has $\alpha[\Lambda(\varphi)]=0$, \emph{i.e.} that $\Psi(\alpha)(\varphi)=0$. We may now conclude that $\supp\Psi(\alpha)\subseteq K$.
In order to obtain statement (ii), fix $v\in \calD( U ,\bbC^n)$ satisfying $\|v\|_{\infty}\leq 1$ and compute:
\begin{eqnarray*}
\|\epsilon \Lambda({\diver_{\bsL^*} v})\|_K&=&\epsilon \|\Lambda({\diver_{\bsL^*}v})\|_K,\\&=&\epsilon\sup\left\{\left|\int_ U  {g}\,{\diver_{\bsL^*} v}\right|:g\in BV_{\bsL,K},\|D_\bsL g\|\leq 1\right\},\\
&=&\epsilon\sup\left\{\left|\int_U \bar{v}\cdot d[D_\bsL {g}]\right|: g\in BV_{\bsL, K}, \|D_\bsL g\|\leq 1\right\},\\
&\leq &\epsilon \sup\{\|D_\bsL g\|\cdot \|v\|_\infty:g\in BV_{\bsL, K}, \|D_\bsL g\|\leq 1\},\\
&\leq &\epsilon,
\end{eqnarray*}
so that one has $\epsilon\Lambda({\diver_{\bsL^*}v})\in V(K,\epsilon)$. It hence follows that:
$$
\epsilon |\Psi(\alpha)(\diver_{\bsL^*} v)|=|\alpha[\epsilon\Lambda({\diver_{\bsL^*} v})]|\leq 1,
$$
and we thus get:
$$
|\Psi(\alpha)(\diver_{\bsL^*} v)|\leq\frac{1}{\epsilon}.
$$
Since $v\in \calD( U ,\bbC^n)$ is an arbitrary vector field satisfying $\|v\|_\infty\leq 1$, this yields $\|D_\bsL \Psi(\alpha)\|\leq\frac{1}{\epsilon}$, and concludes the proof of the claim.
\end{proof}

We now turn to proving Lemma~\ref{lem.inv}.
\begin{proof}[Proof of Lemma~\ref{lem.inv}]
To prove part (i), fix $g\in BV_{\bsL,c}( U )$ and compute, for $\varphi\in\calD( U )$:
$$
\Psi[\Phi(g)](\varphi):=\Phi(g)[\Lambda({\varphi})]=\Lambda({\varphi})({g})=\int_ U  g\bar{\varphi},
$$
that is, $\Psi[\Phi(g)]=g$ in the sense of distributions.

In order to prove part (ii), fix $\alpha\in CH_{\bsL}^*( U )$. We have to show that, for any $F\in CH_\bsL( U )$, we have:
$$
\Phi[\Psi(\alpha)](F)=\alpha(F),
$$
\emph{i.e.} that for any $F\in CH_\bsL( U )$, one has:
$$
F[\Psi(\alpha)]=\alpha(F).
$$
To this purpose, define for any $F\in CH_\bsL( U )$ a map:
$$
\Delta_F:CH_\bsL( U )^*\to\bbC, \alpha\mapsto \Delta_F(\alpha):=F[\Psi(\alpha)].
$$

\begin{Claim}\label{cl.bg}
Given $F\in CH_\bsL( U )$, the map $\Delta_F$ is weakly$^*$-continuous on $V(K,\epsilon)^\circ$ for all $K\subset\subset U $ and $\epsilon>0$.
\end{Claim}
To prove this claim, fix $K\subset\subset  U $, $\epsilon>0$ and assume that $(\alpha_i)_{i\in I}\subseteq$ is a net weak$^*$-converging to $0$. In particular one gets:
\begin{itemize}
\item[(a)] for any $\varphi\in\calD( U )$, we have $\Lambda(\varphi)\in CH_\bsL( U )$ and hence the net $(\Psi(\alpha_i)(\varphi))_{i\in I}=(\alpha_i[\Lambda(\varphi)])_{i\in I}$ converges to $0$.
\end{itemize}
According to Claim~\ref{cl.w*}, we moreover have:
\begin{itemize}
\item[(b)] $\supp\Psi(\alpha_i)\subseteq K$ for each $i\in I$;
\item[(c)] $c:=\sup_{i\in I} \|D_\bsL \Psi(\alpha_i)\|\leq\frac{1}{\epsilon}$.
\end{itemize}
{It hence follow from Proposition~\ref{prop.compacite} that the net $(\|\Psi(\alpha_i)\|_{L^1})_{i\in I}$ converges to $0$}. From the fact that $F$ is an $\bsL$-charge we see that the net $(F[\Psi(\alpha_i)])_{i\in I}$ converges to $0$  as well. This means, in turn, that $(\Delta_F(\alpha_i))_{i\in I}$ converges to $0$, which shows that $\Delta_F$ is weak$^*$-continuous on $V(K,\epsilon)$.

\begin{Claim}
For any $\alpha\in CH_\bsL( U )^*$, we have $\Delta_F(\alpha)=\alpha(F)$.
\end{Claim}
To prove the latter claim, observe that according to Claim~\ref{cl.bg} and to the Banach-Grothendieck theorem \cite[Theorem~8.5.1]{EDW}, there exists $\tilde{F}\in CH_\bsL( U )$ such that for any $\alpha\in CH_\bsL( U )^*$, we have:
$$
\Delta_F(\alpha)=\alpha(\tilde{F}).
$$
Yet given $g\in BV_{\bsL,c}( U )$, we then have, according to [Lemma~\ref{lem.inv}, (i)]:
$$
F(g)=F\{\Psi[\Phi(g)]\}=\Delta_F[\Phi(g)]=\Phi(g)(\tilde{F})=\tilde{F}(g),
$$
\emph{i.e.} $F=\tilde{F}$, which proves the claim.

It now suffices to observe that Lemma~\ref{lem.inv} is proven for we have established the equality $F[\Psi(\alpha)]=\alpha(F)$ for any $F\in CH_\bsL( U )$ and $\alpha\in CH_\bsL( U )^*$.
\end{proof}

As a corollary, we get a proof of the density of $\Gamma[C( U ,\bbC^n)]$ in $CH_\bsL( U )$.
\begin{Corollary}
The space $\Lambda[\calD( U )]$ is dense in $CH_\bsL( U )$.
\end{Corollary}
\begin{proof}
Assuming that $\alpha\in CH_\bsL( U )^*$ satisfies $\alpha\upharpoonright \Lambda[\calD( U )]=0$, we compute for any $\varphi\in\calD( U )$:
$$
\Psi(\alpha)(\varphi):=\alpha[\Lambda(\varphi)]=0.
$$
This means that $\Psi(\alpha)=0$, and implies that $\alpha=\Phi\circ\Psi (\alpha)=\Phi(0)=0$. The result then follows from the Hahn-Banach theorem.
\end{proof}
\begin{Corollary}\label{cor.dens-2}
The space $\Gamma[C( U ,\bbC^n)]$ is dense in $CH_\bsL( U )$.
\end{Corollary}
\begin{proof}
It follows from the previous corollary that $\Lambda[\calD( U )]$ is dense in $CH_\bsL( U )$. Since by hypothesis we also have $\Lambda[\calD(U)]\subseteq \Gamma[C(U,\bbC^n)]\subseteq CH_\bsL(U)$, it is clear that $\Gamma(U,\bbC^n)$ is dense in $CH_\bsL(U)$.
\end{proof}

In order to study the range of $\Gamma^*$, we introduce the following linear operator:
$$
\Xi:BV_{\bsL,c}( U )\to C( U ,\bbC^n)^*, g\mapsto \Xi(g),
$$
defined by $\Xi(g)(v):=\Gamma(v)(g)$ for any $v\in C( U ,\bbC^n)$.

\begin{Claim}
We have $\im\Gamma^*=\im \Xi$.
\end{Claim}
\begin{proof}
To prove this claim, fix $\mu\in C( U ,\bbC^n)$. If one has $\mu=\Gamma^*(\alpha)$ for some $\alpha\in CH_\bsL( U )^*$, then we compute for $v\in C( U ,\bbC^n)$:
$$
\Xi[\Psi(\alpha)](v)=\Gamma(v)[\Psi(\alpha)]=\Phi[\Psi(\alpha)][\Gamma(v)]=\alpha[\Gamma(v)]=\Gamma^*(\alpha)(v)=\mu(v),
$$
so that one has $\mu=\Xi[\Psi(\alpha)]\in \im\Xi$. Conversely, if one has $\mu=\Xi(g)$ for some $g\in BV_{\bsL,c}( U )$, then we compute for $v\in C( U ,\bbC^n)$:
$$
\Gamma^*[\Phi(g)](v)=\Phi(g)[\Gamma(v)]=\Gamma(v)(g)=\Xi(g)(v)=\mu(v),
$$
so that one has $\mu=\Gamma^*[\Phi(g)]\in\im\Gamma^*$.
\end{proof}

Consider the set
$$
B:=\{v\in C( U ,\bbC^n): \|v\|_\infty\leq 1\}.
$$
It is clear that $B$ is bounded in $C( U ,\bbC^n)$. Hence the seminorm:
$$
p:C( U ,\bbC^n)^*\to \R_+, \mu\mapsto p(\mu):=\sup_{v\in B} |\mu(v)|,
$$
is strongly continuous (\emph{i.e.} continuous with respect to the strong topology) on $C( U ,\bbC^n)^*$. Observe now that one has, for $g\in BV_{\bsL,c}( U )$:
\begin{eqnarray*}
p[\Xi(g)]&=&\sup_{v\in B} |\Xi(g)(v)|,\\
&=&\sup\{|\Gamma(v)(g)|:v\in B\},\\
&=& \|D_\bsL g\|.
\end{eqnarray*}
\begin{Lemma}\label{lem.closed}
The set $\im \Xi$ is strongly sequentially closed in $C( U ,\bbC^n)^*$.
\end{Lemma}
\begin{proof}
Fix a sequence $(\Xi(g_k))_{k\in\N}\subseteq\im \Xi$ and assume that, in the strong topology, one has:
$$
\Xi(g_k)\to \mu\in C( U ,\bbC^n)^*,\quad k\to\infty.
$$
The strong continuity of $p$ then yields:
$$
c:=\sup_{k\in\N} \|D_\bsL g_k\|=\sup_{k\in\N} p[\Xi(g_i)]<+\infty.
$$

\begin{Claim}
There exists a compact set $K\subset\subset U $ such that one has $\supp g_k\subseteq K$ for each $k\in\N$.
\end{Claim}
To prove this claim, let us first prove that the sequence $(\supp D_\bsL g_k)_{k\in\N}$ is compactly supported in $ U $ (\emph{i.e.} that there is a compact subset of $ U $ containing $\supp Dg_k$ for all $k$). To this purpose, we proceed towards a contradiction and assume that it is not the case. Let then $ U =\bigcup_{j\in\N}  U _j$ be an exhaustion of $ U $ by open sets satisfying, for each $j\in\N$, $\bar{ U }_j\subseteq  U _{j+1}$ and such that $\bar{ U }_j$ is a compact subset of $ U $ for each $j\in\N$. Since $(\supp D_\bsL g_k)_{k\in\N}$ is not compactly supported, there exist increasing sequences of integers $(j_l)_{l\in\N}$ and $(k_l)_{l\in\N}$ satisfying, for any $l\in\N$:
$$
\supp (D_\bsL g_{k_l})\cap ( U _{j_l+1}\setminus\bar{ U }_{j_l})\neq \emptyset.
$$
In particular, there exists for each $l\in\N$ a vector field $v_l\in C_c( U_{j_{l}+1}\setminus\bar{U}_{j_l} ,\bbC^n)$ with $\|v_l\|_\infty\leq 1$ and:
$$
a_l:=\left|\int_ U  \bar{v}_l\cdot d[Dg_{k_l}]\right|>0.
$$
Let now, for $l\in\N$, $b_l:=\max_{0\leq k\leq l} \frac{1}{a_k}$ and define a bounded set $B'\subseteq C( U ,\bbC^n)$ by:
$$
B':=\left\{v\in C( U ,\bbC^n): \|v\|_{\infty,\bar{ U }_{j_l+1}} \leq l b_l\text{ for each }l\in\N\right\}.
$$
It follows from the construction of $B$ that one has $w_l:=lb_lv_l\in B$ for any $l\in\N$. Moreover the seminorm
$$
p':=C( U ,\bbC^n)^*\to\R_+, \mu\mapsto \sup_{v\in B'} |\mu(v)|,
$$
is strongly continuous. Yet we get for $l\in\N$:
$$
p'[\Xi(g_{k_l})]\geq |\Xi(g_{k_l})(w_l)|=|\Gamma(w_l)(g_{k_l})|=lb_l \left|\int_{ U } \bar{v}_l\cdot d(Dg_{k_l})\right|=lb_la_l\geq l.
$$
Since this yields $p'[\Xi(g_{k_l})]\to\infty$, $l\to\infty$, we get a contradiction with the fact that $p'$ is strongly continuous (recall that $(\Xi(g_{k_l}))_{l\in\N}$ converges in the strong topology).

Now fix $k\in\N$ and $x\in U\setminus\supp (D_\bsL g_k)$; choose an open set $V\subseteq \Omega$ such that one has $V\cap\supp (D_\bsL g_k)=\emptyset$ and observe that one has $\|g_k\|_{L^{N/N-1}(V)}\leq \|D_\bsL g_k\|(V)=0$. It hence follows that $g_k$ is a.e. equal to $0$ on $V$, and hence that $x\notin\supp g_k$. This proves the inclusion $\supp g_k\subseteq K$ for all $k$, which establishes the claim.\\

Getting back to the proof of Lemma~\ref{lem.closed}, observe that, according to Proposition~\ref{prop.compacite}, there exists a subsequence $(g_{k_l})\subseteq (g_k)$, $L^1$-converging to $g\in BV_{\bsL,c}( U )$.
Using the fact that $\Gamma(v)$ is an $\bsL$-charge, we compute:
$$
\mu(v)=\lim_{l\to\infty} \Xi(g_{k_l})(v)=\lim_{l\to\infty} \Gamma(v)(g_{k_l})=\Gamma(v)(g)=\Xi(g)(v),
$$
and hence we get $\mu=\Xi(g)\in\im \Xi$.
\end{proof}

We hence proved the following theorem.
\begin{Theorem}
We have $CH_\bsL( U )=\Gamma[C( U ,\bbC^n)]$.
\end{Theorem}

\section{Appendix}

\begin{theorem}\label{l1}
Let $p(x,D)$ be a pseudodifferential operator with symbol in the H\"ormander class $S^{m}_{1,0}(\R^{N})$ and consider $k(x,y)$ be the distribution kernel of $p(x,D)$ defined by the oscillatory integral
\begin{equation}
k(x,y)=\int e^{{2i\pi}(x-y)\cdot \xi}p(x,\xi)d\xi.
\end{equation}
If $m<0$ then $p(x,D)$ maps continuously $L^{1}(\R^{N})$ onto itself. 
\end{theorem}

\begin{proof}
Writing $p(x,D)u=(k(x,\cdot)\ast u)(x)$ it is sufficient to prove that $k(x,y) \in L^{1}(\R^{N}\times \R^{N})$ using a pointwise control of the kernel due to \`Alvarez and Hounie in \cite{AH}. In order to prove the boundedness in $\LL^{1}$ norm, we first localize the kernel in the diagonal region. Let  $A=\left\{(x,y)\in \R^{N}\times \R^{N} \; : \;|x-y|< 1\right\}$ {be a} neighborhood of the diagonal. % and complementary %$\left\{|x-y|\geq 1\right\}$ .  
% Now we consider $(x,y) \in A$. 
If $m<-N$ then $k$ is bounded and clearly   the property follows. If $0<m+N$ then there {exists} $C>0$ such that $|k(x,y)|\leq C |x-y|^{-(m+N)}$, and then $k$ is integrable {on $A$}, since $m<0$. The limiting case occurs when $m=-N${, which implies} $|k(x,y)|\leq C \log|x-y|$ {from which the property follows}. On the other hand,  by the pseudo-local property (see  \cite[Theorem 1.1]{AH}), {we see that} there exists $L_{0} \in \Z^{+}$ such that $|k(x,y)|\leq |x-y|^{-L}$ for $L\geq L_{0}$ and $|x-y| \geq 1${; hence} $k$ is integrable{ on $\complement A$}, since $L\geq \max \left\{N,L_{0} \right\}$. Combining {all those} cases we conclude that $k(x,y) \in L^{1}(\R^{N}\times \R^{N})$.
\end{proof}

Consider a class of pseudodifferential operators, called Bessel potential $J_{\beta}$ for {$\beta >0$},  defined by 
$$J_{\beta}f(x)=\int_{\R^{N}}e^{{2i\pi}  x \cdot \xi}b(x,\xi)\hat{f}(\xi)d\xi, \quad f \in S'(\mathbb{R}^{N}),$$
where $b(x,\xi)=\langle \xi \rangle^{\beta}:=(1+4\pi^{2}|\xi|^{2})^{-\beta/2}$ belongs to the H\"ormander class $S^{-\beta}_{1,0}(\R^{N})$. We define the nonhomogeneous Sobolev space $W^{\beta,p}(\R^{N})$ for $\beta>0$ and $1\leq p<\infty$ as
$$W^{\beta,p}(\R^{N})=\left\{ f \in S'(\R^{N}) : J_{-\beta}f \in L^{p}(\R^{N}) \right\}$$
with associated norm $\|f\|_{{k,p}}:=\|J_{-\beta}f\|_{p}$.
As a consequence of {Theorem 3.5 in} \cite{AH} and Theorem \ref{l1} when $p=1$, it follows that $W^{\beta,p}(\R^{N}) \subset L^{p}(\R^{N})$ continuously, \emph{i.e.} that one has:
\begin{equation}\label{ineq0.1}
\|u\|_{{p}}=\|J_{\beta}(J_{-\beta}u)\|_{{p}}\leq C \|J_{-\beta}u\|_{{p}}=C \|u\|_{{\beta,p}}.
\end{equation}
For $B=B(x_0,\ell)$ a fixed ball let $\tilde{B}=B(x_0,2\ell)$ the ball with the same center as $B$ but twice its radius. Let $\psi \in C_{c}^{\infty}(\tilde{B})$ satisfy $\psi(x) \equiv 1$ on $B$ and define  $\Lambda_{\beta} := \Lambda_{\beta}(x,D)$ the pseudodifferential operator with symbol $\lambda_{\beta}(x,\xi)=\psi(x)\left\langle \xi\right\rangle^{\beta}$. 
Denote by $W^{\beta,p}_{c}(B)$ the set of distributions $f \in \mathcal{E}'(B)$ such that $\Lambda_{\beta}f \in \LL^{p}(\R^{N})$, endowed with the semi-norm $\|f\|_{{\beta,p}(B)}\;:=\;\|\Lambda_{\beta}u\|_{p}$. Note that the space  $W^{\beta,p}_{c}(B)$ is independent of the choice of $\psi$. %; i.e, if $\psi_{2}(x), \psi_{1}(x) \in C_{c}^{\infty}(\widetilde{B})$ satisfies $\psi_{1}(x)=\psi_{2}(x) \equiv 1$ on $B$ then  $\|\Lambda_{\beta, \psi_{1}}f\|_{L^{p}}\; \cong \;\ \|\Lambda_{\beta, \psi_{2}}f\|_{L^{p}}$.  
 In view of \eqref{ineq0.1}, we have the continuous inclusion:
$$ W^{\beta,p}_{c}(B) \subset \LL^{p}(\R^{N}),$$%L^{p}_{c}{(\tilde{B})}:=L^{1}(\R^{N}) \cap \mathcal{E}'(\tilde{B}).$$ 
for $1\leq p<\infty$. Next we present a version of the Rellich-Kondrachov compactness for $W_{c}^{\beta,1}(B)$.
\begin{theorem}\label{thm.compact-beta}
Let $0<\beta<1$. The embedding $W_{c}^{\beta,1}(B)\subset \subset \LL^{1}(\R^{N})$ is compact.
\end{theorem}
The proof follows the same strategy as the proof of Theorem~A in \cite{HKP} and will be presented for the sake of completeness. The compact embedding of $W_{c}^{\beta,p}(B)$ in $\LL^p(\R^n)$ for $1<p<\infty$ could be established by analogous means.    

\begin{proof}

According to the previous comments on continuity, it is enough to verify the compactness.  We will show that if $(u_m)$ is a bounded sequence in {$W_{c}^{\beta,1}(B)$} then there exist a subsequence $(u_{m_j})_j$ which converges in {$\LL^{1}(\R^{N})$}. 
Consider the regularizations $u^\epsilon_m=\eta_\epsilon*u_m$ where $\eta\in C_{c}^{\infty}(B_0^1)$, {$\int_{\R^N}\eta=1$}, $\eta_\epsilon(x)=\epsilon^{-N}\eta(x/\epsilon)$ and $0<\epsilon\le1$. It is enough to show 
that the family $\{ u^\epsilon_m\}_{\epsilon,m}$ has the following two properties:
\begin{enumerate}
  \item [(i)] for any fixed $0<\epsilon<1$, the sequence $(u^\epsilon_m)_{m\in\N}$
        is a {relatively compact} subset of $\LL_{c}^1(B'):=\LL^{1}(\R^{N}) \cap \mathcal{E}'(B')$;
  \item [(ii)] $u^\epsilon_m \to u_m$ in $\LL_{c}^{1}(B')$ {\it uniformly} in $m$
                   as $\epsilon \searrow 0$,
\end{enumerate}  
where $B'$ is a closed ball that contains the support of all $u^{\epsilon}_{m}$.  

Since the inclusion $C_{c}(B') \subset \LL^{1}_{c}(B')$ is continuous, property (i) will follow once we shall have proven that {$(u^{\epsilon}_{m})_m$} is a precompact subset of $C_{c}(B')$. We claim that for each $\epsilon>0$, {$(u^{\eps}_{m})_m$} is uniformly bounded and equicontinuous. In fact, one has for $x\in B'$:
\begin{eqnarray*}
|u^{\eps}_{m}(x)|&=&\left|\left\langle u_{m}, \eta_{\eps}(x-\cdot)\right\rangle\right|,\\
&\leq &\|\Lambda_{\beta}\Lambda_{-\beta}u_{m}\|_{1}\|\eta_{\eps}\|_{\infty},\\ 
&\leq &C(B)\|\ \Lambda_{-\beta} u_{m}\|_{\LL^{1}(B)}\|\eta_{\eps}\|_{\infty},\\
&\leq &C({B}) \eps^{-N}\|u_{m}\|_{{\beta,1}},
\end{eqnarray*}
and analogously
\begin{eqnarray*}
|\nabla u^{\eps}_{m}(x)| &\leq& \| \Lambda_{-\beta}u_{m}\|_{{1}}\|(\Lambda_{\beta}\circ\nabla) \eta_{\eps}\|_{\infty}\\
&\leq &C(B) \eps^{-(N+1-\beta)}\|u_{m}\|_{{\beta,1}}.
\end{eqnarray*}
The conclusion follows from Arzel\`a-Ascoli theorem. 

To prove (ii)  we will first consider the identity :
\begin{eqnarray*}\nonumber
 u^{\eps}_{m}(x)-u_{m}(x)&=&\int_{0}^{\eps}\frac{\partial}{\partial s}\left(u_{m} \ast \eta_{s}\right)(x)\,ds \notag ,\\
&=&- \int_{0}^{\eps}\left\{u_{m} \ast {\nabla \cdot [x\eta]_{s}}\right\}(x)\,ds \notag,\\
&=&- \int_{0}^{\eps}\left\{\Lambda_{-\beta}u_{m} \ast { (\Lambda_{\beta} \circ \nabla) \cdot [x\eta]_{s}}\right\}(x)\,ds.
\end{eqnarray*}
But from the equalities {$\Gamma_{\beta}(t,\xi):={2i\pi }\psi(x)\sum_{k=1}^{N}\xi_{k}(t^{2}+4\pi^{2}|\xi|^{2})^{-\frac{\beta}{2}}$} we get:
$${(\Lambda_{\beta}\circ \nabla)\cdot g_{s}(x) }  % \int_{\R^{N}} e^{2\pi i x \cdot \xi}\, \xi^{\alpha}\;\widehat{g}(s\xi)\,d \xi
=s^{{\beta-1}}[\Gamma_{\beta}(s,D)g]_{s}(x),$$
after which we compute, using Fubini's theorem:
\begin{eqnarray*}
\int_{\R^{N}}\big| u^{\eps}_{m} - u_{m}\big|(x) &\leq& \int_{\R^{N}} \int_{0}^{\eps}s^{{\beta-1}}\left| \Lambda_{-\beta}u_{m} \ast (\Gamma_{\beta}(s,D) [{x\eta}])_{s}\right|(x)\,dsdx,\\
&\leq&  \int_{\R^{N}}\int_{0}^{\eps}s^{{\beta-1}} \left(\int_{K \subset \R^{n}}|(\Gamma_{\beta}(s,D){[y\eta]})_{{s}}(y)|\cdot|\Lambda_{-\beta}u_{m}(x-{y})|dy\right)dsdx,\\
&\leq&  \int_{0}^{\eps} \int_{K \subset \R^{N}}s^{{\beta-1}}|(\Gamma_{\beta}(s,D) [{y\eta}])_{{s}}(y)|\cdot \left(\int_{\R^{N}}|\Lambda_{-\beta}u_{m}(x-{y})|dx\right)dyds,\\
& \leq & C\eps^{\gamma}\| \Lambda_{-\beta}u_{m}\|_{{1}}.
\end{eqnarray*}
{To obtain the latter inequalities, we observe (defining $\tilde{\Gamma}_\alpha(t,\xi):={2i\pi}  \sum_{k=1}^{N}\xi_{k} (t^2+4\pi^2|\xi|^2)^{-\frac{\alpha}{2}}$ and letting $\tilde{B}$ be the ball defined above):}
\begin{align*}
\bigg| \int_{0}^{\eps}\int_{K\subset\R^{N}} \!\!s^{\alpha-1}[\Gamma_{\alpha}( s,D)g]_{{s}}(y)dyds \bigg| &= \left|  \int_{0}^{\eps} \left( \int_{{\tilde{B}}} \psi\left(\frac{{y}}{s} \right) \left[{ \tilde{\Gamma}_{\alpha}(1,D)g_{s}} \right] (y)dy \right)ds \right|, \\	
%& =  \left|  \int_{0}^{\eps} \left( \int_{\tilde{K}\subset\R^{N}} \left[ \Gamma_{\alpha}(1,D)g_{s} \right] (y)dy \right)ds \right| \\
	%& \leq  c(\tilde{K})   \int_{0}^{\eps} \left\| \Gamma_{\alpha}(1,D)g_{s}\right\|_{L^{r}}  ds \\
	& \leq C  \int_{0}^{\eps} \left\| g_{s} \right\|_{{r}} ds,  \\
	%& \lesssim   \left\| g \right\|_{L^{r}} \int_{0}^{\eps} s^{\frac{N}{r}-N} ds  \\
	& \leq C'  \eps^\gamma,
\end{align*}
{where $C=C(\eta,\tilde{K})>0$ and $\gamma=\frac{N}{r}-N+\alpha>0$ are constants for $1<r<N/(N-\alpha)$.}

To finish the proof, we claim that, for a given $\delta>0$, there exists a subsequence {$(u_{m_j})_j \subset (u_{m})_m$} such that one has:
\begin{equation}\label{b3}
\limsupe_{j,k \rightarrow \infty} \|u_{m_{j}}-u_{m_k}\|_{1}\leq \delta.
\end{equation}
Indeed, for $\eps>0$ sufficiently small,  we have:
\begin{equation}\label{b1}
\|u^{\eps}_{m}-u_{m}\|_{1}\leq \delta/2
\end{equation}
uniformly in $m$. Since $(u_{m})$ and $(u^{\eps}_{m})$
are supported in a closed ball $B'$, by Arzel\`a-Ascoli's theorem there exists a subsequence 
$(u^{\eps}_{m_j})_j$ wich converges uniformly in $B'$. In particular, this yields:
\begin{equation}\label{b2}
\limsupe_{j,k \rightarrow \infty} \|u^{\eps}_{m_{j}}-u^{\eps}_{m_k}\|_{{1}}=0.
\end{equation}
Note that \eqref{b3} is a consequence of \eqref{b1} and \eqref{b2}. 
Using \eqref{b3} for $\delta=1/n$ for $n=1,2,3,...$ and the diagonal process  we can extract a convergent subsequence {$(u_{m_\ell})_\ell$}.

\end{proof}

\section*{Acknowledgments}
We wish to thank Prof. Jorge Hounie (Universidade Federal de S\~ao Carlos) for helpful suggestions concerning this work.

\bigskip

{\small
\parbox[t]{3.5in}
{\textbf{Laurent Moonens}\\
Universit\'e Paris-Sud\\
Laboratoire de Math\'ematique UMR~8628\\
Universit\'e Paris-Saclay\\
B\^atiment 425\\
F-91405 Orsay Cedex (France)\\
E-mail: Laurent.Moonens@math.u-psud.fr}
\bigskip

{\small
\parbox[t]{5in}
{\textbf{Tiago Picon}\\
University of S\~ao Paulo\\
Faculdade de Filosofia, Ci\^encias e Letras de Ribeir\~ao Preto\\ 
Departamento de Computa\c{c}\~ao e Matem\'atica\\
Avenida Bandeirantes 3900, CEP 1404-040,
Ribeir\~ao Preto, Brasil\\
E-mail: picon@ffclrp.usp.br}
\bigskip

\end{document}